\documentclass[12pt,twoside]{article}

\newtheorem{definition}{Definition}
\newtheorem{theorem}{Theorem}
\newtheorem{lemma}{Lemma}
\newtheorem{corollary}{Corollary}
\newtheorem{remark}{Remark}

\newtheorem{proposition}{Proposition}
\usepackage{amsmath}
\usepackage{amssymb}
\usepackage{enumitem}

\begin{document}
\date{}
\title{
\vspace{60px}
\large\textbf{Generalizations of the Numerical Radius, Crawford Number and Numerical Index Functions in the Weighted Case}\\
}

\author{\small{{\bf Zameddin I. Ismailov$^{1}$, Pembe Ipek Al$^{1}$ }
\let\thefootnote\relax\footnote{{$^{1}$Karadeniz Technical University, Department of Mathematics, 61080, Trabzon, Turkey.
\newline E-mail : \texttt{zameddin.ismailov@gmail.com,}
\texttt{ipekpembe@gmail.com} }}}}

\date{}

\maketitle

\vspace{0.6cm}

\noindent
\small{\textbf{
\begin{center}
Abstract
\end{center}
}}
In this article, firstly, some simple and smoothness properties of the weighted numerical radius and the weighted Crawford number functions are investigated. Then, some generalization formulas for lower and upper bounds of the weighted numerical radius function are obtained. Later on, some evaluations for lower and upper bounds of the weighted numerical index are given. The obtained results are generalized some well-known famous results about the special weighted numerical radius and the special weighted Crawford number functions in the recently literature. Also, the different and useful results are provided to the literature. 
\vspace{0.125cm}
\footnotesize
\noindent
\hspace{0,38cm}

\begin{tabbing}
$\mathit{Key Words:}$ \textup {weighted numerical radius, weighted Crawford number, weighted numerical index.} \vspace*{0.5cm} 
\end{tabbing}

\noindent

\section{Introduction}
Throughout this article, $ H $ denotes a complex Hilbert space endowed with the inner product $ \langle \cdot , \cdot \rangle $ and associated norm $ \Vert \cdot \Vert . $ Let $ \mathbb{B} (H) $ stand for the $ C^{*} $-algebra of all bounded linear operators acting on $ H. $  The numerical radius of an operator $ A $ is given by
$$
\omega (A) =\sup \lbrace \vert\langle Ax,x\rangle\vert: \ x\in H, \ \Vert x \Vert=1 \rbrace .
$$

The usual operator norm and the Crawford number of an operator $ A $  are, respectively, defined by
$$
\Vert A \Vert=\sup \lbrace \Vert Ax \Vert: \ x\in H, \  \Vert x \Vert=1 \rbrace
$$
and 
$$
c(A)=\inf \lbrace \vert\langle Ax,x\rangle\vert:x\in H, \ \Vert x \Vert=1 \rbrace 
$$
\cite{Halmos}.

Remember that for any $ A\in \mathbb{B}(H) $ the classical numerical radius $ \omega (A) $ is a norm on $ \mathbb{B}(H) $ and equivalent to the operator norm that satisfies the relation 
\begin{eqnarray}
\label{equ1}
\frac{1}{2} \Vert A \Vert \leq \omega (A) \leq \Vert A \Vert
\end{eqnarray}
(see \cite{Halmos}).

For any $ A\in \mathbb{B}(H) $ one of improvement formula for the numerical radius in form 
\begin{eqnarray}
\label{equ2}
\frac{1}{4}\Vert A^{*}A+AA^{*}\Vert \leq \omega ^{2}(A)\leq \frac{1}{2} \Vert A^{*}A+AA^{*}\Vert
\end{eqnarray}
has been obtained in \cite{Kittaneh}. 

Assume that $ \varphi:[0,1]\rightarrow \mathbb{R} $ and $ \psi:[0,1]\rightarrow \mathbb{R} $ are continuous functions. For the $ A\in \mathbb{B}(H) $ and $ t\in [0,1] $ the weighted numerical radius and the weighted Crawford number functions will be defined by 
$$
\omega_{t}(\varphi,\psi;A)=\sup\limits _{x\in S_{1}(H)}\vert \langle \left( \varphi(t)A+\psi(t)A^{*}\right)x,x  \rangle  \vert
$$
and
$$
c_{t}(\varphi,\psi;A)=\inf\limits _{x\in S_{1}(H)}\vert \langle \left( \varphi(t)A+\psi(t)A^{*}\right)x,x  \rangle  \vert,
$$
respectively, where $ S_{1}(H) $ is a unite sphere in $ H. $

Similarly, the weighted operator norm of any operator $ A\in \mathbb{B}(H) $ will be defined as 
$$
\Vert A \Vert_{t}=\Vert \varphi(t)A+\psi(t)A^{*}\Vert ,t\in [0,1].
$$

It is clear that if $ \varphi (t)=1, \ \psi(t)=0, \ t\in [0,1] $, then $ \omega_{t}(1,0;A) $ and $ c_{t}(1,0;A) $  coincide with the classical numerical radius $ \omega (A) $ and the classical Crawford number $ c(A) $, respectively.

In case when $ \varphi (t)=1, \ \psi(t)=1-2t, \ 0\leq t\leq 1 $ and $ \varphi (t)=\nu, \ \psi(t)=1-\nu, \ 0\leq \nu\leq 1, $  these special weighted numerical radii have been investigated in \cite{Nayak} and \cite{Sheikhhosseini}, respectively. There in some upper bounds for the weighted numerical radius have been researched too. 

If we take $ \varphi, \ \psi \in C[0,1] $ real-valued functions,
$$
A=\begin{pmatrix}
        0 & 1\\
        3 & 0
    \end{pmatrix}, \ A:\mathbb{C}^{2}\rightarrow \mathbb{C}^{2},
$$
then
$$
\varphi A+ \psi A^{*}=\begin{pmatrix}
        0 & \varphi +3 \psi\\
        3\varphi +\psi & 0
    \end{pmatrix}.
$$
For $ x=(x_{1},x_{2})\in \mathbb{C}^{2}, $ we have
$$
\langle(\varphi A+\psi A^{*})x,x\rangle=(\varphi + 3\psi)x_{2}\overline{x_{1}} + (3\varphi + \psi)x_{1}\overline{x_{2}}.
$$
Hence, we get
$$
\vert(\varphi A+\psi A^{*})x\vert^{2}=\vert(\varphi + 3\psi)x_{2}\vert^{2} + \vert(3\varphi + \psi)x_{1}\vert^{2}.
$$
In this case, from the result Lemma 2 in \cite{Abu} for any $ t\in [0,1] $  we have
$$
\omega_{t}(\varphi,\psi;A)=\frac{1}{2}\sup\limits_{\theta \in \mathbb{R}}\vert e^{i\theta} (\varphi +3\psi )+e^{-i\theta}  (3\varphi +\psi ) \vert = \frac{1}{2}\sup\limits_{\theta \in \mathbb{R}}\vert \varphi (e^{i\theta}+3e^{-i\theta})+ \psi(3e^{i\theta}+e^{-i\theta}) \vert .
$$
Also,
$$
c_{t}(\varphi,\psi;A)=0 \ \ \text{and} \ \ \Vert A \Vert_{t}=\vert (\varphi+3\psi )(3\varphi +\psi ) \vert (t).
$$

For more improvements results related with lower and upper bounds of the classical numerical radius we refer to studies in recent years (see, \cite{Baklouti}, \cite{Bhunia1}, \cite{Bhunia2}, \cite{Feki1}, \cite{Feki2}, \cite{Kittaneh}, \cite{Zamani}). 

In mathematical literature, also there is a constant on a Banach space, known as the numerical index of the space, which relates the behaviour of the numerical radius with the usual norm of an operator. The numerical index of the Banach space $ X $ is the constant
$$
n(X)= \inf \left\lbrace \omega (A): A\in B(X), \Vert A \Vert=1 \right\rbrace ,
$$
equivalent, $ n(X) $ is the maximum of those $ k\geq 0 $ such that $ k \Vert A \Vert \leq \omega (A) $ for every $ A\in B(X). $ This notion was introduced and studied in the 1970 paper \cite{Duncan}, see also the monographs \cite{Bonsall1}, \cite{Bonsall2} and the survey paper \cite{Kadets} for background. Clearly, $ 0\leq n(X)\leq 1, \ n(X)>0 $ means that the numerical radius is a norm on $ B(X) $ equivalent to the operator norm and $ n(X)=1 $ if and only if numerical radius and operator norm coincide. If $ X $ is a complex Banach space, then $ \frac{1}{e}\leq n(X)\leq 1 $ and if $ X $ is a real Banach space, then $ 0\leq n(H)\leq 1. $ If $ H $ is a complex Hilbert space, then $ n(H)=\frac{1}{2} $ and for a real Hilbert space $ H, \ n(H)=0 .$ Moreover, $ n(l_{1})=n(l_{1}^{m})=n(l_{\infty} )=n(l_{\infty}^{m})=1, $ where $ m\in \mathbb{N}. $ All these results can be found in \cite{Bonsall1} and \cite{Kadets}. Some recent developments for the study of the numerical index are \cite{Choi}, \cite{Martín2},  \cite{Martín1}, \cite{Merí}, \cite{Sain}.

The main aim of this study is to generalize some well-known results about the weighted numerical radius and the weighted Crawford number functions  in the mathematical literature \cite{Bhunia3}, \cite{Halmos}, \cite{Kittaneh}, \cite{Nayak},  \cite{Sheikhhosseini}, as well as to provide different and useful results to the this area. 

This work is organized as follows: In Section 2, some simple and smoothness properties of the weighted numerical radius and the weighted Crawford number functions have been investigated. In Section 3, some evolutions formulas for lower and upper bounds of the weighted numerical radius function have been obtained. Later on, in Section 4, some evaluations for lower and upper bounds of the weighted numerical index have been given. These results are generalization of some well-known results in the literature.

Note that each operator $ A\in \mathbb{B}(H) $ can be expressed in the Cartesian decomposition form as $ A=ReA+iImA, $ where $ ReA=\frac{A+A^{*}}{2} $ and $ ImA=\frac{A-A^{*}}{2i}. $ Here, $ A^{*} $ denote the adjoint of $ A. $ Throughout this study it will be denote by $ \vert A \vert= (A^{*}A)^{1/2} $ the absolute value of an operator $ A\in \mathbb{B}(H) . $

\section{Some properties of the weighted numerical radius and the weighted Crawford number functions}

Let us begin this section with some simple properties of the weighted numerical radius and the weighted Crawford number functions.

\begin{proposition}
\label{pro2.1}
For any $ A, B \in \mathbb{B}(H) $ and $ t\in [0,1] $ the followings are true:
\begin{enumerate}[label=(\arabic*)]
\item $ \omega_{t}(0,\psi;A)=\vert \psi(t) \vert \omega(A) $ and  $ c_{t}(0,\psi;A)=\vert \psi(t) \vert c(A), $
\item $ \omega_{t}(\varphi,0;A)=\vert \varphi(t) \vert \omega(A) $ and  $c_{t}(\varphi,0;A)=\vert \varphi(t) \vert c(A),  $
\item $ \omega_{t}(1,1;A)=2\Vert ReA \Vert $ and $ c_{t}(1,1;A)=2\Vert ReA \Vert, $
\item $ \omega_{t}(1,-1;A)=2\Vert ImA \Vert $ and $ c_{t}(1,-1;A)=2\Vert ImA \Vert, $
\item $ \omega_{t}(\varphi,\psi;iA)=\omega_{t}(\varphi,-\psi;A), $ 
\item $ \frac{1}{2}\Vert A \Vert_{t}\leq \omega_{t}(\varphi,\psi;A)\leq  \Vert A \Vert_{t}, $
\item $ \omega_{t}(\varphi,\psi;A)\leq \left( \vert \varphi \vert + \vert \psi \vert\right)(t)\omega (A) $ and $ c_{t}(\varphi,\psi;A) \geq \vert \vert \varphi \vert - \vert  \psi \vert \vert (t) c (A),$
\item $ \omega_{t}(\varphi,\varphi;A)=2 \vert \varphi(t)\vert \omega(ReA)=2\vert \varphi(t)\vert \Vert ReA \Vert, $
\item $ \omega_{t}(\psi,\psi;A)=2 \vert \psi(t) \vert\omega (ReA)=2 \vert \psi(t) \vert\Vert ReA \Vert, $
\item If $ A=A^{*}, $ then $ \omega_{t}(\varphi,\psi;A)= \vert  \varphi + \psi \vert (t) \omega (A)=  \vert  \varphi + \psi \vert (t) \Vert A \Vert , $
\item $ \omega_{t}(\varphi,\psi;A)=\omega_{t}(\psi,\varphi;A^{*}), $ 
\item If $ A=A^{*}, $ then $ \omega_{t}(\varphi,\psi;A)=\omega_{t}(\psi, \varphi;A),$
\item $ \omega_{t}(\varphi,\psi;A+B)\leq\omega_{t}(\varphi,\psi;A)+\omega_{t}(\varphi,\psi;B),$
\item $ c_{t}(\varphi,\psi;A+B)\leq\omega_{t}(\varphi,\psi;A)+c_{t}(\varphi,\psi;B),$
\item $ \omega_{t}(\varphi,\psi;AB)\leq \left( \vert \varphi \vert + \vert \psi \vert \right) (t)\omega(AB),$ 
\item $ c_{t}(\varphi,\psi;AB)\leq\vert \varphi \vert (t) c(AB) + \vert \psi \vert (t)\omega (AB),$
\item $ \Vert \varphi AB+\psi (AB)^{*} \Vert^{2}_{t}\leq \left( \vert \varphi \vert^{2} +  \vert \psi \vert^{2}\right)(t)\Vert AB \Vert^{2}+\vert \varphi \psi \vert (t)\omega ((AB)^{2}+  \vert \varphi \psi \vert (t)\omega ( (B^{*}A^{*})^{2}, $
\item $ \omega_{t}(\varphi_{1}+\varphi_{2}, \psi_{1}+\psi_{2};A)\leq \omega_{t}(\varphi_{1},\psi_{1};A)+\omega_{t}(\varphi_{2},\psi_{2};A),$
\item $ c_{t}(\varphi,\psi;A)\geq  \inf \left\lbrace \vert \varphi +\psi \vert(t), \vert \varphi -\psi \vert(t) \right\rbrace m(A), $ 
where 
$
m(A)=\inf\limits_{x\in S_{1}(H)}\vert (Ax,x) \vert .
$
\end{enumerate}
\end{proposition}
\noindent {\it Proof.}
To give an idea it will be proved the 15th, 17th and 19th claims. Firstly, let us start with the proof of the 15th claim. For any $ A, B \in \mathbb{B}(H) $ and $ t\in [0,1], $ we get
\begin{eqnarray}
\vert \langle (\varphi (t)AB+\psi (t)B^{*}A^{*})x,x \rangle \vert^{2} & = & \vert \langle \varphi (t)ABx,x \rangle + \langle \psi (t)B^{*}A^{*}x,x \rangle \vert^{2} \nonumber \\
& = & \vert \varphi \vert^{2}(t) \vert \langle ABx,x \rangle\vert^{2} +(\varphi\psi)(t) \left[ \langle ABx,x\rangle^{2}+ \langle x,ABx\rangle^{2} \right]  + \vert \psi \vert^{2}(t) \vert \langle x,ABx \rangle\vert^{2} \nonumber \\
& \leq & \left( \vert \varphi \vert^{2}(t)  +2\vert \varphi\psi\vert(t) + \vert \psi \vert^{2}(t)\right)  \vert \langle ABx,x \rangle\vert^{2} \nonumber \\
& = & \left( \vert  \varphi \vert + \vert\psi \vert\right) ^{2}(t) \vert \langle ABx,x \rangle\vert^{2}. \nonumber
\end{eqnarray}
Then,
$$
\vert \langle (\varphi (t)AB+\psi (t)B^{*}A^{*})x,x \rangle \vert\leq \left( \vert  \varphi \vert + \vert\psi \vert\right) (t) \vert \langle ABx,x \rangle\vert .
$$
Hence, for any $ t\in [0,1] $ we have
$$
\omega_{t}(\varphi,\psi;AB)=  \left( \vert  \varphi \vert + \vert\psi \vert\right) (t)   \omega(AB).
$$
Now, let prove the 17th claim. For $ t\in [0,1], \ x\in S_{1}(H) $ and $ A\in \mathbb{B}(H) $, we have
\begin{eqnarray}
\Vert (\varphi (t)AB + \psi (t) B^{*}A^{*})x \Vert_{t}^{2} & = & \vert \varphi^{2} \vert (t)\langle ABx, ABx\rangle + (\varphi \psi )(t) \langle ABx, B^{*}A^{*}x\rangle +  (\varphi \psi )(t) \langle B^{*}A^{*}x, ABx\rangle \nonumber \\
&  & + \vert \psi \vert^{2}(t)  \langle B^{*}A^{*}x, B^{*}A^{*}x\rangle  \nonumber \\
& \leq & \vert \varphi^{2} \vert (t)\Vert ABx\Vert^{2} + \vert\varphi \psi \vert(t) \omega ( (AB)^{2})+  \vert\varphi \psi \vert(t)  \omega ( (B^{*}A^{*})^{2})
+ \vert \psi \vert^{2}(t)  \Vert B^{*}A^{*}x\Vert^{2}.\nonumber
\end{eqnarray}
Thus, we get
$$
\Vert \varphi AB+\psi (AB)^{*} \Vert^{2}_{t}\leq \left( \vert \varphi \vert^{2} +  \vert \psi \vert^{2}\right)(t)\Vert AB \Vert^{2}+\vert \varphi \psi \vert (t)\omega ( (AB)^{2})+  \vert \varphi \psi \vert (t)\omega ( (B^{*}A^{*})^{2}).
$$
If we take $ \varphi (t)=1 $ and $ \psi (t)=1-2t, \ 0\leq t\leq 1 $ in the last result, we get the inequality proved in [26, Th.
 2.6]
$$
\Vert (1-2t)(AB)^{*} + AB \Vert^{2} \leq (2-4t+4t^{2}) \Vert AB \Vert^{2} +(1-2t) \omega ( (AB)^{2}) + (1-2t)\omega ( (B^{*}A^{*})^{2}). 
$$
So, the property (17) in Proposition \ref{pro2.1} generalizes [26, Th. 2.6].

Lastly, let prove the 19th claim. For $ t\in [0,1], \ x\in S_{1}(H) $ and $ A\in \mathbb{B}(H) $, with the following simple calculations we have
\begin{eqnarray}
\vert \langle(\varphi A+ \psi A^{*})x,x\rangle\vert^{2}  & = &  \vert \langle(\varphi ReA+ \psi ReA)x,x\rangle+i(\varphi ImA - \psi ImA)x,x)\vert^{2} \nonumber \\
& = & \vert \varphi + \psi \vert^{2}(t)\vert \langle ReAx,x\rangle \vert^{2}+\vert \varphi - \psi \vert^{2}(t)\vert \langle ImAx,x\rangle \vert^{2}. \nonumber 
\end{eqnarray}
Then, we get
$$ 
c_{t}^{2}(\varphi,\psi;A)\geq \left( \inf \left\lbrace \vert \varphi +\psi \vert(t), \vert \varphi -\psi \vert(t) \right\rbrace  \right)^{2}m^{2}(A).
$$
\begin{lemma}
\label{lem2.2}
Let $ A,B\in \mathbb{B}(H) $ and $ t\in [0,1], $ then
\begin{enumerate}[label=(\arabic*)]
\item $ \vert \omega_{t}(\varphi,\psi;A)-\omega_{t}(\varphi,\psi;B) \vert \leq \omega_{t}(\varphi,\psi;A-B) $,
\item $ \vert c_{t}(\varphi,\psi;A)-c_{t}(\varphi,\psi;B) \vert \leq c_{t}(\varphi,\psi;A-B) $.
\end{enumerate}
\end{lemma}
\noindent {\it Proof.}
Using the property (13) of Proposition \ref{pro2.1}, we have 
$$
\omega_{t}(\varphi,\psi;A)\leq \omega_{t}(\varphi,\psi;A-B)+\omega_{t}(\varphi,\psi;B) \ \ \text{and} \ \ 
\omega_{t}(\varphi,\psi;B)\leq \omega_{t}(\varphi,\psi;A-B)+\omega_{t}(\varphi,\psi;A)\nonumber
$$
for any $ A, B\in \mathbb{B}(H)$ and $ t\in [0,1]. $ Then, the last inequalities imply that 
$$
\vert \omega_{t}(\varphi,\psi;A)-\omega_{t}(\varphi,\psi;B) \vert \leq \omega_{t}(\varphi,\psi;A-B) .
$$

Similarly, using the property (14) of Proposition \ref{pro2.1}, we have
$$
c_{t}(\varphi,\psi;A)\leq \omega_{t}(\varphi,\psi;A-B)+c_{t}(\varphi,\psi;B) \ \ \text{and} \ \
c_{t}(\varphi,\psi;B)\leq \omega_{t}(\varphi,\psi;A-B)+c_{t}(\varphi,\psi;A)
$$
for any $ A, B\in \mathbb{B}(H) $ and $ t\in [0,1]. $ Then, the last inequalities imply that 
$$
\vert c_{t}(\varphi,\psi;A)-c_{t}(\varphi,\psi;B) \vert \leq \omega_{t}(\varphi,\psi;A-B) .
$$
for any $ A, B\in \mathbb{B}(H) $ and $ t\in [0,1].$

Now, let give smoothness properties of the weighted numerical radius and the weighted Crawford number functions.
\begin{theorem}
\label{thm2.3}
If the function sequences $ (\varphi_{n}) $ and $ (\psi_{n}) $ pointwise converge to functions $ \varphi:[0,1]\rightarrow \mathbb{R} $ and $ \psi:[0,1]\rightarrow \mathbb{R}, $ respectively, then for any $ A\in \mathbb{B}(H) $ and $ t\in [0,1]$
$$
\omega_{t}(\varphi,\psi;A)=\lim\limits_{n\rightarrow \infty}\omega_{t}(\varphi_{n},\psi_{n};A) \ \  \text{and } \ \
c_{t}(\varphi,\psi;A)=\lim\limits_{n\rightarrow \infty}c_{t}(\varphi_{n},\psi_{n};A).
$$
\end{theorem}
\noindent {\it Proof.}
Since 
\begin{eqnarray}
\vert \langle \left( \varphi(t)A+\psi(t)A^{*}\right)x,x  \rangle \vert & \leq & \vert \langle \left( (\varphi-\varphi_{n})(t)A+(\psi-\psi_{n})(t)A^{*}\right)x,x  \rangle  \vert + \vert \langle \left( \varphi_{n}(t)A+\psi_{n}(t)A^{*}\right)x,x  \rangle  \vert \nonumber
\end{eqnarray}
and
\begin{eqnarray}
\vert \langle \left( \varphi_{n}(t)A+\psi_{n}(t)A^{*}\right)x,x  \rangle  \vert & \leq & \vert \langle \left( (\varphi-\varphi_{n})(t)A+(\psi-\psi_{n})(t)A^{*}\right)x,x  \rangle \vert + \vert \langle \left( \varphi(t)A+\psi(t)A^{*}\right)x,x  \rangle  \vert \nonumber
\end{eqnarray}
for any $ A\in \mathbb{B}(H), \ t\in [0,1] $, $ x\in S_{1}(H) $ and $ n\geq 1, $ then 
$$
\vert\omega_{t}(\varphi_{n},\psi_{n};A)-\omega_{t}(\varphi,\psi;A)\vert\leq \omega_{t}(\varphi_{n}-\varphi,\psi_{n}-\psi;A)\leq \left( \vert \varphi_{n}-\varphi \vert(t)+ \vert \psi_{n}-\psi \vert(t) \right) \Vert A \Vert
$$
and 
$$
\vert c_{t}(\varphi_{n},\psi_{n};A)-c_{t}(\varphi,\psi;A)\vert\leq \omega_{t}(\varphi_{n}-\varphi,\psi_{n}-\psi;A)\leq \left( \vert \varphi_{n}-\varphi \vert(t)+ \vert \psi_{n}-\psi \vert(t) \right) \Vert A \Vert.
$$
From the last inequalities and convergence of sequences $ (\varphi_{n}) $ to $ \varphi $ and $ (\psi_{n}) $ to $ \psi, $ the validity of theorem is clear.
\begin{definition}\cite{Kato}
A sequence $ (A_{n}) \subset \mathbb{B}(H) $ is said to  uniformly converges to $ A\in \mathbb{B}(H), $ if for any $ \epsilon >0 , $ there exists a positive integer $ N $ such that for all $ n\in \mathbb{N} $ 
$$
\Vert A_{n}- A \Vert <\epsilon .
$$
\end{definition}
\begin{theorem}
\label{thm2.5}
If the operator sequences $ (A_{n}) $ in $ \mathbb{B}(H) $ uniformly converges on norm $ \Vert \cdot \Vert_{t}, \ 0\leq t \leq 1 $ to the operator $ A\in \mathbb{B}(H), $ then 
$$
\omega_{t}(\varphi,\psi;A)=\lim\limits_{n\rightarrow \infty}\omega_{t}(\varphi,\psi;A_{n}) \ \ \text{and} \ \
c_{t}(\varphi,\psi;A)=\lim\limits_{n\rightarrow \infty}c_{t}(\varphi,\psi;A_{n}).
$$
\end{theorem}
\noindent {\it Proof.}
For $ t\in [0,1] $ and $ n\geq 1 $ by property (1) of Lemma \ref{lem2.2}, we have
\begin{eqnarray}
\vert \omega_{t}(\varphi, \psi; A_{n})-\omega_{t}(\varphi, \psi; A) \vert \leq  \omega_{t}(\varphi, \psi; A_{n}-A)
 \leq  \Vert \varphi(t)(A_{n}-A)+\psi(t)(A_{n}-A)  \Vert 
 = \Vert A_{n}-A\Vert_{t}, \nonumber
\end{eqnarray}
and similarly,
\begin{eqnarray}
\vert c_{t}(\varphi, \psi; A_{n})-c_{t}(\varphi, \psi; A) \vert  \leq  \Vert A_{n}-A\Vert_{t}. \nonumber
\end{eqnarray}
Then, the claims of theorem are clear.
\begin{theorem}
\label{thm2.6}
If $ \varphi, \psi \in H_{\alpha}[0,1], \ 0<\alpha \leq 1, $ then 
$
\omega_{t}(\varphi, \psi; A), \ c_{t}(\varphi, \psi; A)\in H_{\alpha}[0,1]
$ for any $ A\in \mathbb{B}(H) $,
where $ H_{\alpha}[0,1] $ is class of Hölder functions with degree $ \alpha \in (0,1] $ in $ [0,1] $.
\end{theorem}
\noindent {\it Proof.}
For any $ 0\leq t,s \leq 1, $ we have 
\begin{eqnarray}
\omega_{t}(\varphi, \psi; A) & = & \omega (\varphi(t)A+\psi(t)A^{*}) \nonumber \\
& = & \omega (\varphi(s)A+\psi(s)A^{*})+((\varphi(t)-\varphi(s))A+(\psi(t)-\psi(s))A^{*}) \nonumber \\
& \leq & \omega (\varphi(s)A+\psi(s)A^{*})+ \vert \varphi(t)-\varphi(s)\vert + \vert \psi(t)-\psi(s)\vert \Vert A \Vert \nonumber \\
& \leq & \omega_{s}(\varphi, \psi; A)+K_{\varphi}\vert t-s\vert ^{\alpha} + K_{\psi}\vert t-s\vert ^{\alpha} \nonumber \\
& \leq & \omega_{s}(\varphi, \psi; A)+K\vert t-s\vert ^{\alpha}, \ K=\max\lbrace K_{\varphi}, K_{\psi} \rbrace  \nonumber 
\end{eqnarray}
and similarly,
$$
\omega_{s}(\varphi, \psi; A)\leq \omega_{t}(\varphi, \psi; A)+ K\vert t-s\vert ^{\alpha}.
$$
Consequently, for each $ 0\leq t,s\leq 1 $ and $ A\in \mathbb{B}(H) $ the following inequality is true
$$
\vert \omega_{t}(\varphi, \psi; A)-\omega_{s}(\varphi, \psi; A) \vert\leq K\vert t-s\vert ^{\alpha},
$$
i.e.
$$
 \omega_{t}(\varphi, \psi; A)\in H_{\alpha}[0,1].
$$
By similarly way the validity of second claim can be proved.
\begin{theorem}
\label{thm2.7}
Let $ A\in \mathbb{B}(H), $ then $  \omega_{t}(\varphi, \psi; A), \  c_{t}(\varphi, \psi; A)\in C[0,1], $ where $ C[0,1] $ is class of continuous functions on $ [0,1]. $
\end{theorem}
\noindent {\it Proof.}
For $ A\in \mathbb{B}(H) $ and $ t,s \in [0,1] $  we have 
\begin{eqnarray}
\vert \langle  (\varphi (t)A+\psi (t)A^{*})x,x \rangle   \vert & = & \vert \langle  ((\varphi (t)-\varphi (s))A+(\psi (t)-\psi (s))A^{*})x,x \rangle + \langle  (\varphi (s)A+\psi (s)A^{*})x,x \rangle  \vert \nonumber \\
& \leq & \vert \langle  ((\varphi (t)-\varphi (s))A+(\psi (t)-\psi (s))A^{*})x,x \rangle \vert + \vert \langle  (\varphi (s)A+\psi (s)A^{*})x,x \rangle  \vert. \nonumber 
\end{eqnarray}
Let $ s $ be any fixed number in $ [0,1]. $ Hence, we get 
$$
\omega_{t}(\varphi,\psi;A)\leq \omega(\varphi(t)-\varphi(s),\psi(t)-\psi(s);A^{*})+ \omega_{s}(\varphi,\psi;A).
$$
Similarly,
$$
\omega_{s}(\varphi,\psi;A)\leq \omega(\varphi(t)-\varphi(s),\psi(t)-\psi(s);A^{*})+ \omega_{t}(\varphi,\psi;A).
$$
Consequently, from the last two relations and property (7) of Proposition \ref{pro2.1}
we get
\begin{eqnarray}
\vert \omega_{t}(\varphi, \psi; A)-\omega_{s}(\varphi, \psi; A) \vert & = & \vert  \omega (\varphi(t)A+\psi(t)A^{*})-\omega(\varphi(s)A+\psi(s)A^{*}) \vert \nonumber \\
& \leq & \omega ((\varphi(t)-\varphi(s))A+(\psi(t)-\psi(s))A^{*}) \nonumber \\
& \leq & \left( \vert \varphi(t)-\varphi(s)\vert + \vert \psi (t)- \psi(s) \vert \right)  \omega (A) \nonumber   \\
& \leq & \left( \vert \varphi(t)-\varphi(s)\vert + \vert \psi (t)- \psi(s) \vert \right)  \Vert A \Vert . \nonumber  
\end{eqnarray}
Since the functions $ \varphi $ and $ \psi $ is continuous on $ [0,1] $, then the continuity of $ \omega_{t}(\varphi, \psi; A), \ 0\leq t \leq 1 $ is clear by the last relation.

Similarly, from the subadditivity of the classical numerical radius function the continuity of function $ c_{t}(\varphi, \psi; A) $ on $ [0,1] $ can be easily proved.

Later on, from the definitions of the weighted numerical radius and the weighted Crawford numbers function it implies the validity of next result.

\begin{proposition}
\label{pro2.8}
If $ \varphi, \ \psi \in D[0,1] $ and $ \varphi\geq 0, \ \psi \geq 0, $ then for every $ A\in \mathbb{B}(H) $
\begin{enumerate}[label=(\arabic*)]
\item $ \omega_{t}'(\varphi,0;A)=\omega_{t}(\varphi',0;A)$,
\item $ \omega_{t}'(0,\psi;A)=\omega_{t}(0,\psi';A)$,
\item $ c_{t}'(\varphi,0;A)=c_{t}(\varphi',0;A)$,
\item $ c_{t}'(0,\psi;A)=c_{t}(0,\psi';A)$,
\end{enumerate}
where $ D[0,1] $ is class of differentiable functions on $ [0,1] $.
\end{proposition}

\begin{proposition}
\label{pro2.9}
For $ A\in \mathbb{B}(H) $, the followings are true:
\begin{enumerate} [label=(\arabic*)]
\item If $ \varphi \neq -\psi , $ then 
$$
\Vert Re A \Vert \leq \left( \int\limits_{0}^{1}\vert \varphi + \psi \vert (t) dt \right)^{-1} \int\limits_{0}^{1} \omega_{t}(\varphi,\psi;A) dt ,
$$
\item If $ \varphi \neq \psi , $ then 
$$
\Vert Im A \Vert \leq \left( \int\limits_{0}^{1}\vert \varphi - \psi \vert (t) dt \right)^{-1} \int\limits_{0}^{1} \omega_{t}(\varphi,\psi;A) dt ,
$$
\item If $ \vert \varphi \vert \neq \vert \psi \vert $, then
$$
\Vert A \Vert \leq \left[ \left( \int\limits_{0}^{1}\vert \varphi + \psi \vert (t) dt \right)^{-1}+\left( \int\limits_{0}^{1}\vert \varphi - \psi \vert (t) dt \right)^{-1} \right] \int\limits_{0}^{1} \omega_{t}(\varphi,\psi;A) dt . 
$$
\end{enumerate}
\end{proposition}
\noindent {\it Proof.}
For any $ x\in S_{1}(H) $ , we have 
\begin{eqnarray}
\vert \varphi + \psi \vert (t) \vert \langle Re Ax,x\rangle \vert   \leq  \sqrt{\vert \varphi + \psi \vert^{2} (t)\vert \langle Re Ax,x\rangle \vert^{2}+\vert \varphi - \psi \vert^{2} (t)\vert \langle Im Ax,x\rangle \vert^{2}} 
 =  \vert \langle(\varphi A+ \psi A^{*})x,x\rangle \vert \nonumber 
\end{eqnarray}
and 
\begin{eqnarray}
\vert \varphi - \psi \vert (t) \vert \langle Im Ax,x\rangle \vert   \leq  \sqrt{\vert \varphi + \psi \vert^{2} (t)\vert \langle Re Ax,x\rangle \vert^{2}+\vert \varphi - \psi \vert^{2} (t)\vert \langle Im Ax,x\rangle \vert^{2}} 
 =  \vert \langle(\varphi A+ \psi A^{*})x,x\rangle \vert     . \nonumber 
\end{eqnarray}
Hence, we get
$$
\vert \varphi +\psi \vert (t) \Vert ReA \Vert \leq \omega_{t}(\varphi, \psi, A) \ \ \text{and} \ \
\vert \varphi -\psi \vert (t) \Vert ImA \Vert \leq \omega_{t}(\varphi, \psi, A).
$$
Consequently, from the last inequalities  we have 
\begin{enumerate} [label=(\arabic*)]
\item If $ \varphi \neq -\psi , $ then  
$$
\Vert Re A \Vert \leq \left( \int\limits_{0}^{1}\vert \varphi + \psi \vert (t) dt \right)^{-1} \int\limits_{0}^{1} \omega_{t}(\varphi,\psi;A) dt ,
$$
\item If $ \varphi \neq \psi , $ then 
$$
\Vert Im A \Vert \leq \left( \int\limits_{0}^{1}\vert \varphi - \psi \vert (t) dt \right)^{-1} \int\limits_{0}^{1} \omega_{t}(\varphi,\psi;A) dt .
$$
\item From the first and second claims of this theorem, it is clear that the third claim is true.
\end{enumerate}

Using the result obtained for the classical numerical radius and the classical Crawford number in \cite{Otkun}, the following theorem is obtained for the weighted numerical radius and the weighted Crawford number functions.

\begin{theorem}
If for any $ n\geq 1, \ H_{n}$ is a Hilbert space, $ A_{n}\in \mathbb{B}(H_{n}) $, $ H=\bigoplus\limits_{n=1}^{\infty}H_{n} $ and $ A=\bigoplus\limits_{n=1}^{\infty}A_{n}, \ A\in \mathbb{B}(H), $ then 
\begin{enumerate} [label=(\arabic*)]
\item $ \omega_{t}(\varphi,\psi;A)=\sup\limits_{n\geq 1} \omega_{t}(\varphi,\psi;A_{n})$,
\item In case when $ Re(A_{n})\geq 0 $ (or $ Re(A_{n})\leq 0 $), $ n\geq 1 $, $ c_{t}(\varphi,\psi;A)=\inf\limits_{n\geq 1} \omega_{t}(\varphi,\psi;A_{n})$.
\end{enumerate}
\end{theorem}

\section{On the lower and upper bounds of the weighted numerical radius}
In this section, some estimates for lower and upper bounds of the weighted numerical radius are given. Firstly, we give some well-known auxiliary results of work \cite{Bhunia3}.

\begin{lemma}
\label{lem 3.1}
Let $ A\in \mathbb{B}(H), $ then
\begin{enumerate} [label=(\arabic*)]
\item $\omega (A)\geq \frac{\Vert A \Vert}{2} + \frac{\vert \ \Vert ReA \Vert-\Vert Im A \Vert \ \vert }{2}, $
\item $\omega^{2} (A)\geq \frac{1}{4}\Vert A^{*}A+AA^{*}\Vert + \frac{\vert \ \Vert ReA \Vert^{2}-\Vert Im A \Vert^{2} \ \vert }{2}, $ 
\item $\omega^{2} (A)\geq \frac{1}{4}\Vert A^{*}A+AA^{*}\Vert + \frac{c^{2} (ReA) + c^{2} (ImA)}{2} + \left\vert \frac{ \Vert ReA \Vert^{2}-\Vert Im A \Vert^{2} }{2}+ \frac{c^{2} (ImA) - c^{2} (ReA)}{2} \right\vert , $
\item  $\omega^{4} (A)\geq \frac{1}{16}\Vert (A^{*}A+AA^{*})^{2}+4Re^{2}(A^{2})\Vert + \frac{1}{2}\left\vert \Vert ReA \Vert^{4}-\Vert ImA \Vert^{4}  \right\vert . $
\end{enumerate}
\end{lemma}

We obtain the following lower bound for the weighted numerical radius of bounded linear operators.
\begin{theorem}
\label{thm 3.2}
Let $ A\in \mathbb{B}(H)$ and $ t\in [0,1] $, then
\begin{enumerate} [label=(\arabic*)]
\item $\omega_{t} (\varphi,\psi;A)\geq \inf\left\lbrace \vert \varphi + \psi \vert (t), \vert \varphi - \psi \vert (t)  \right\rbrace \left( 
\frac{\Vert A \Vert}{2} + \frac{\vert \ \Vert ReA \Vert-\Vert Im A \Vert \ \vert }{2}\right) , $
\item $\omega_{t}^{2} (\varphi,\psi;A)\geq \left( \inf\left\lbrace \vert \varphi + \psi \vert (t), \vert \varphi - \psi \vert (t)  \right\rbrace\right) ^{2} \left( \frac{1}{4}\Vert A^{*}A+AA^{*}\Vert + \frac{\vert \ \Vert ReA \Vert^{2}-\Vert Im A \Vert^{2} \ \vert }{2}\right) , $ 
\item 
\begin{eqnarray}
\omega_{t}^{2} (\varphi,\psi;A) & \geq & \left( \inf\left\lbrace \vert \varphi + \psi \vert (t), \vert \varphi - \psi \vert (t)  \right\rbrace\right) ^{2} \left(  \frac{1}{4}\Vert A^{*}A+AA^{*}\Vert + \frac{c^{2} (ReA) + c^{2} (ImA)}{2}\right.  \nonumber \\
&  & + \left. \left\vert \frac{ \Vert ReA \Vert^{2}-\Vert Im A \Vert^{2} }{2}+ \frac{c^{2} (ImA) - c^{2} (ReA)}{2} \right\vert \right) \nonumber ,
\end{eqnarray}
\item  $\omega_{t}^{4} (\varphi,\psi;A)\geq \left( \inf\left\lbrace \vert \varphi + \psi \vert (t), \vert \varphi - \psi \vert (t)  \right\rbrace \right) ^{4} \left( \frac{1}{16}\Vert (A^{*}A+AA^{*})^{2}+4Re^{2}(A^{2})\Vert + \frac{1}{2}\left\vert \Vert ReA \Vert^{4}-\Vert ImA \Vert^{4}  \right\vert \right) . $
\end{enumerate}
\end{theorem}
\noindent {\it Proof.}
(1) For any $ x\in S_{1}(H), $ we have
\begin{eqnarray}
\vert \langle (\varphi A+\psi A^{*})x,x\rangle\vert^{2} & = & \vert \langle ((\varphi + \psi ) ReA+(\varphi-\psi) ImA)x,x\rangle\vert^{2} \nonumber \\
& = & \vert \varphi + \psi \vert^{2} \vert \langle ReAx,x\rangle \vert^{2}+ \vert \varphi - \psi \vert^{2} \vert \langle ImAx,x\rangle \vert^{2} \nonumber \\
& \geq & \left( \inf\left\lbrace \vert \varphi + \psi \vert (t), \vert \varphi - \psi \vert (t)  \right\rbrace\right) ^{2} \left( \vert \langle ReAx,x\rangle \vert^{2}+  \vert \langle ImAx,x\rangle \vert^{2} \  \right). \nonumber
\end{eqnarray}
Then, 
$$
\omega_{t} (\varphi, \psi; A)\geq \inf\left\lbrace \vert \varphi + \psi \vert (t), \vert \varphi - \psi \vert (t)  \right\rbrace  \max\left\lbrace \Vert ReA \Vert, \Vert ImA \Vert\right\rbrace. 
$$
From the last inequality and property (1) of Lemma \ref{lem 3.1}, the validity of first claim is clear.

(2) For any $ x\in S_{1}(H) $, we get
$$
\vert \langle (\varphi A+\psi A^{*})x,x\rangle\vert^{2} \geq \left( \inf\left\lbrace \vert \varphi + \psi \vert (t), \vert \varphi - \psi \vert (t)  \right\rbrace\right) ^{2} \left( \vert \langle ReAx,x\rangle \vert^{2}+  \vert \langle ImAx,x\rangle \vert^{2} \  \right).
$$

Thus, we have 
$$
\omega_{t}^{2} (\varphi,\psi;A)\geq \left( \inf\left\lbrace \vert \varphi + \psi \vert (t), \vert \varphi - \psi \vert (t)  \right\rbrace \right)  ^{2} \max\left\lbrace \Vert ReA \Vert^{2}, \Vert ImA \Vert^{2}\right\rbrace. 
$$
Then, from the last inequality and property (2) of Lemma \ref{lem 3.1} the validity of second claim is clear.

(3) For any $ x\in S_{1}(H), $ we have
\begin{eqnarray}
\vert \langle (\varphi A+\psi A^{*})x,x\rangle\vert^{2} & = & \vert \langle \vert(\varphi + \psi \vert^{2} \vert \langle ReAx,x \rangle\vert^{2} + \vert \varphi -\psi \vert^{2}\vert \langle ImAx,x\rangle \vert^{2}\nonumber \\
& \geq & \left( \inf\left\lbrace \vert \varphi + \psi \vert (t), \vert \varphi - \psi \vert (t)  \right\rbrace\right) ^{2} \max \left\lbrace    
\Vert Re A \Vert^{2}+c^{2}(ImA), \Vert Im A \Vert^{2}+c^{2}(ReA)
\right\rbrace  \nonumber.
\end{eqnarray}
From the last inequality, definition of the weighted numerical radius and property (3) of Lemma \ref{lem 3.1} it is obtained third claim of theorem.

(4) For any $ x\in S_{1}(H) $ it implies that
$$
\vert \langle (\varphi A+\psi A^{*})x,x\rangle\vert^{4} \geq \left( \inf \left\lbrace \vert \varphi + \psi \vert (t), \vert \varphi - \psi \vert (t)  \right\rbrace\right) ^{4} \max \left\lbrace \Vert Re A \Vert^{4}, \Vert Im A \Vert^{4} \right\rbrace .
$$
From the last inequality and property (4) of Lemma \ref{lem 3.1} it is obtained fourth claim of theorem.

\begin{corollary}
If we take $ \varphi (t)=1, \ \psi (t)=0, \ 0\leq t \leq 1 $ in Theorem \ref{thm 3.2}, then Theorem \ref{thm 3.2} and Lemma \ref{lem 3.1} coincide. So, Theorem \ref{thm 3.2} generalizes Lemma \ref{lem 3.1}.
\end{corollary}

\noindent {\it Example.} In the complex Hilbert space $ L^{2}(0,1) $ consider the following classical Volterra integration operator 
$$
V: L^{2}(0,1)\rightarrow L^{2}(0,1), \ Vf(x)= \int\limits_{0}^{x}f(t)dt, \ f\in L^{2}(0,1) .
$$
It is known that $ \Vert V \Vert =\frac{2}{\pi}, \ \Vert ReV \Vert =\frac{1}{2}, $ and $ \Vert ImV \Vert =\frac{1}{\pi} $ (see \cite{Khadkhuu1}). Then, by Theorem \ref{thm 3.2}, we have 
$$\omega_{t} (\varphi,\psi;A)\geq \inf\left\lbrace \vert \varphi + \psi \vert (t), \vert \varphi - \psi \vert (t)  \right\rbrace \frac{2+\pi}{4\pi}.
$$
Noted that $ \omega (V)=\frac{1}{2}, \ \omega (ReV)= \frac{1}{2}, $ and $ \omega(ImV)=\frac{1}{\pi} $ (see \cite{Khadkhuu2}).
It is easily shown that
$$
\Vert V^{*}V+VV^{*} \Vert=\frac{2}{\pi^{2}}. 
$$
The operator $ V^{*}V+VV^{*}: L^{2}(0,1)\rightarrow L^{2}(0,1) $ is compact selfadjoint and positive. Then,
$$
\Vert V^{*}V+VV^{*} \Vert = \sup \left\lbrace \lambda : \lambda\in \sigma ( V^{*}V+VV^{*} ) \right\rbrace ,
$$
here $ \sigma(\cdot) $ is defined as the set of spectrum of an operator (see, \cite{Hirsch}). Consider the following spectral problem 
$$
(V^{*}V+VV^{*} )f=\lambda f, \ \lambda\neq 0, \ f\in L^{2}(0,1),
$$
i.e. 
$$
\int\limits_{x}^{1}\int\limits_{0}^{y}f(t)dtdy+\int\limits_{0}^{x}\int\limits_{y}^{1}f(t)dtdy=\lambda f.
$$
If $ f $ is chosen to be a twice differentiable function on (0,1), from the last equation we have
$$
\begin{cases}
\lambda f^{''}=-2f, \\
\lambda (f(0)+f(1))=\int\limits_{0}^{1}f(t)dt, \\
f^{'}(0)+f^{'}(1)=0.
\end{cases}
$$
Assume that
$$
f(t)=\cos\left( \sqrt{\frac{2}{\lambda}}t \right), \ 0\leq t\leq 1. 
$$
Then,
$
f^{''}=\frac{-2}{\lambda}f
$
is satisfied and from the boundary condition $ f^{'}(0)+f^{'}(1)=0 $  we have
$$
\lambda_{n}=\frac{2}{n^{2}\pi^{2}}, \ n\in \mathbb{N}.
$$
In this case $ \max\limits_{n\geq 1}\lambda_{n}=\lambda_{1}=\frac{2}{\pi^{2}}. $

On the other hand,
$$
\int\limits_{0}^{1}f_{1}(t)dt=\int\limits_{0}^{1} cos\left( \sqrt{\frac{2}{\lambda_{1}}}t \right)dt=0, 
$$
$$
f_{1}(0)+f_{1}(1)=1+cos\pi =0.
$$
Then, for $ \lambda_{1}=\frac{2}{\pi^{2}} $ and $ f_{1}(t)=cos\left( \sqrt{\frac{2}{\pi}}t \right), \ 0\leq t \leq 1 $
the condition 
$$
\lambda_{1}(f_{1}(0)+f_{1}(1))=\int\limits_{0}^{1}f_{1}(t)dt
$$
is satisfied. Hence, $ \Vert V^{*}V+VV^{*} \Vert=\frac{2}{\pi^{2}}. $

It is known that $ \Vert V\Vert=\frac{2}{\pi}, \ \Vert ReV \Vert= \frac{1}{2}, $ $ \Vert ImV \Vert=\frac{1}{\pi} $ and $ \omega (V)=\frac{1}{2}, \ c (ReV)= 0, $ and $ c(ImV)=0 $ (see, \cite{Khadkhuu1}, \cite{Khadkhuu2}).
Then, by Theorem \ref{thm 3.2} we have
$$
\omega_{t} (\varphi,\psi;V)\geq \inf\left\lbrace \vert \varphi + \psi \vert (t), \vert \varphi - \psi \vert (t)  \right\rbrace \frac{2+\pi}{4\pi}
$$
and
$$
\omega_{t} ^{2}(\varphi,\psi;V)\geq \frac{1}{8}\inf\left\lbrace \vert \varphi + \psi \vert (t), \vert \varphi - \psi \vert (t)  \right\rbrace .
$$

Now, we give a few well-known inequalities that are essential for proving our theorems, starting with Buzano's inequality.
\begin{lemma} 
\cite{Buzano}
\label{lem 3.4}
 Let $ x,y,e\in H $ with $ \Vert e \Vert=1. $ Then 
$$
\vert \langle x,e\rangle \langle e,y\rangle\vert\leq \frac{1}{2} \left( \vert \langle x,y \rangle \vert+ \Vert x \Vert \Vert y  \Vert\right). 
$$

The next lemma pertains to a positive operator.
\end{lemma}
\begin{lemma} \cite{Nayak}
\label{lem 3.5}
Let $ T\in \mathbb{B}(H) $ be a positive operator. Then for all $ r\geq 1 $ and $ x\in H $ with $ \Vert x \Vert =1,  $ we have 
$$
\langle Tx,x \rangle^{r}\leq \langle T^{r}x,x\rangle .
$$
\end{lemma}

Next, we present the generalized mixed Schwarz inequality.
\begin{lemma}
\label{lem 3.6}
\cite{Furuta} If $ T\in \mathbb{B}(H), $ then for all $ x,y \in H $ and $ \alpha \in [0,1] $
$$
\vert \langle Tx,y\rangle \vert^{2}\leq \langle \vert T \vert^{2\alpha}x,x\rangle \langle \vert T^{*}\vert^{2(1-\alpha )}y,y\rangle .
$$
\end{lemma}

We obtain the following upper bound for the weighted numerical radius of bounded linear operators.
\begin{theorem}
\label{thm 3.7}
Let $ A\in \mathbb{B}(H) $ and $ t\in [0,1], $ then 
$$
\omega_{t}^{2}(\varphi,\psi;A) \leq \vert \varphi \vert^{2}(t)\omega^{2}(A)+\vert \varphi \psi \vert (t)\omega(A^{2})+\frac{(\vert \varphi \vert+\vert \psi \vert)(t)}{2}\vert \psi \vert (t) \Vert A^{*}A+AA^{*}\Vert.
$$
\end{theorem}
\noindent {\it Proof.}
For $ t\in [0,1]  $ and $ x\in S_{1}(H) $, using Lemma \ref{lem 3.4}-\ref{lem 3.6} we have 
\begin{eqnarray}
& & \vert \langle( \varphi (t) A+ \psi (t) A^{*})x,x\rangle\vert^{2} \nonumber \\
& \leq & \vert \varphi \vert ^{2}(t) \vert \langle Ax,x\rangle \vert^{2}+2\vert \varphi \psi \vert (t) \vert \langle Ax,x\rangle\vert \vert \langle A^{*}x,x\rangle\vert + \vert \psi \vert^{2}(t) \vert \langle A^{*}x,x\rangle\vert^{2} \nonumber \\
& = &  \vert \varphi \vert ^{2}(t) \vert \langle Ax,x\rangle \vert^{2}+2\vert \varphi \psi \vert (t)  \vert\langle Ax,x\rangle \langle x,A^{*}x\rangle \vert+ \vert \psi \vert^{2}(t) \vert \langle Ax,x\rangle\vert^{2} \nonumber \\
& \leq & \vert \varphi \vert ^{2}(t) \vert \langle Ax,x\rangle \vert^{2}+\vert \varphi \psi  \vert (t) \left[ \vert \langle Ax,A^{*}x\rangle \vert + \Vert Ax \Vert \Vert A^{*}x \Vert \right] + \vert\psi^{2}\vert (t) \langle \vert A \vert x,x\rangle \langle\vert A^{*}\vert x,x\rangle  \nonumber \\
& \leq & \vert \varphi \vert^{2} (t) \vert \langle Ax,x\rangle\vert^{2}+\vert \varphi \psi \vert (t) \vert \langle A^{2}x,x\rangle \vert + \vert \varphi \psi \vert (t) \frac{1}{2} \langle (\vert A \vert^{2}+ \vert A^{*}\vert^{2})x,x\rangle  + \vert \psi \vert^{2}(t)\frac{1}{2} \langle\vert A \vert^{2}+ \vert A^{*}\vert^{2})x,x\rangle \nonumber \\
& \leq & \vert \varphi \vert^{2}(t)\omega^{2}(A)+\vert \varphi \psi \vert (t)\omega(A^{2})+\frac{(\vert \varphi \vert+\vert \psi \vert)(t)}{2}\vert \psi \vert (t) \Vert A^{*}A+AA^{*}\Vert. \nonumber 
\end{eqnarray}
From the last inequality and the definition of the weighted numerical radius the validity of claim is clear.
\begin{corollary}
\label{cor3.5}
If we take $ \varphi (t)=1-2t, \ \psi (t)=1, \ 0\leq t \leq 1 $, and $ T=A^{*}$ in Theorem \ref{thm 3.7}, we get the inequality proved in [26, Th. 2.4], which says that
\begin{eqnarray}
\omega_{t}^{2}(T)\leq (1-2t)^{2}\omega^{2}(T)+(1-2t)\omega(T^{2})+(1-t)\Vert T^{*}T+TT^{*}\Vert . \nonumber
\end{eqnarray}
So, the inequality obtained in Theorem \ref{thm 3.7} generalizes [26, Th. 2.4].

Also, if we take $ \varphi (t)=0, \ \psi (t)=1, \ 0\leq t \leq 1 $ in Theorem \ref{thm 3.7}, we get the right side of inequality (\ref{equ2}) which says that 
$$
\omega^{2}(A)\leq \frac{1}{2}\Vert A^{*}A+AA^{*} \Vert .
$$
So, the inequality obtained in Theorem \ref{thm 3.7} generalizes the right side of inequality (\ref{equ2}).

\end{corollary}

In the next theorem we obtain another upper bound for the weighted numerical radius.
\begin{theorem}
\label{thm3.9}
Let $ A\in \mathbb{B}(H) $ and $ t\in [0,1], $ then
$$
\omega_{t}^{2}(\varphi ,\psi;A) \leq \frac{1}{2} ( \vert \varphi \vert^{2}+ \vert \psi \vert^{2})(t) \Vert A^{*}A+AA^{*}\Vert + \vert\varphi \psi\vert (t) \omega (A^{2}+(A^{*})^{2}).
$$
\end{theorem}
\noindent {\it Proof.}
For every $ t\in [0,1] $ and $ x\in S_{1}(H) $ using Lemma \ref{lem 3.5} and Lemma \ref{lem 3.6} we have
\begin{eqnarray}
& & \vert \langle(\varphi (t)A+\psi (t) A^{*})x,x\rangle\vert^{2} \nonumber \\
& \leq & \langle \vert \varphi (t)A + \psi (t) A^{*}\vert x,x \rangle \langle \vert \varphi (t)A^{*} + \psi (t) A\vert x,x \rangle \nonumber \\
& \leq & \frac{1}{2}\left[ \langle \vert \varphi (t)A + \psi (t) A^{*}\vert^{2} x,x \rangle + \langle \vert \varphi (t)A^{*} + \psi (t) A\vert^{2} x,x \rangle \right] \nonumber \\
& \leq & \frac{1}{2} \left[  \langle(\varphi (t)A+\psi(t)A^{*})(\varphi(t)A^{*}+\psi(t)A)x,x\rangle + \langle(\varphi (t)A^{*}+\psi(t)A)(\varphi(t)A+\psi(t)A^{*})x,x\rangle\right] \nonumber \\
& = &  \frac{1}{2}  \left[  \langle \vert \varphi \vert^{2}(t)AA^{*}+\varphi (t)\psi (t)AA+\varphi (t)\psi (t) A^{*}A^{*}+ \vert \psi \vert^{2}(t)A^{*}A + \right. \nonumber \\
&  & + \left.  \vert \varphi \vert^{2}(t)A^{*}A+\varphi (t)\psi (t)A^{*}A^{*}+\varphi (t)\psi (t) AA+ \vert \psi \vert^{2}(t)AA^{*})x,x \rangle \right]  \nonumber \\
& = &  \frac{1}{2}  \left[  \langle \vert \varphi \vert^{2}(t)(AA^{*}+A^{*}A)+2\varphi (t)\psi (t)A^{*}A^{*} + 2\varphi (t)\psi (t) AA+ \vert \psi \vert^{2}(t)(A^{*}A+AA^{*}))x,x \rangle\right]   \nonumber \\
& = &  \frac{1}{2}  \left[  \langle(\vert \varphi \vert^{2}(t)+ \vert \psi \vert^{2}(t))(A^{*}A+AA^{*})+2\varphi (t)\psi (t)(A^{*}A^{*}+AA))x,x\rangle \right]   \nonumber 
\end{eqnarray}
From the last relation we have 
$$
\omega_{t}^{2}(\varphi ,\psi;A) \leq \frac{1}{2} ( \vert \varphi \vert^{2}+ \vert \psi \vert^{2})(t) \Vert A^{*}A+AA^{*}\Vert + \vert\varphi \psi\vert (t) \omega (A^{2}+(A^{*})^{2}).
$$
\begin{corollary}
\label{cor3.7}
If we take $ \varphi (t)=1-2t, \ \psi (t)=1, \ 0\leq t \leq 1 $ in Theorem \ref{thm3.9}, we get the inequality proved in [26, Th. 2.7], which says that
$$
\omega_{t}^{2}(T)\leq (1-2t+2t^{2})\Vert T^{*}T+TT^{*}\Vert +(1-2t)\omega(T^{2}+(T^{*})^{2}).
$$
So, the inequality obtained in Theorem \ref{thm3.9} generalizes [26, Th. 2.7].

Also, if we take $ \varphi (t)=0, \ \psi (t)=1, \ 0\leq t \leq 1 $ in Theorem \ref{thm3.9}, we get the right side of inequality (\ref{equ2}).
So, the inequality obtained in Theorem \ref{thm3.9} generalizes the right side of inequality (\ref{equ2}).
\end{corollary}

Now, another estimate from upper bound of the weighted numerical radius will be given.
\begin{theorem}
\label{thm3.11}
Let $ A\in \mathbb{B}(H) $ and $ t\in [0,1], $ then,
$$
\omega^{2}_{t}(\varphi ,\psi;A) \leq \left( \vert \varphi \vert^{2}+ \vert \psi \vert^{2} \right)(t)\omega^{2}(A)+\vert \varphi \psi \vert (t) \omega (A^{2})+ \frac{1}{2}\vert \varphi \psi \vert (t) \Vert AA^{*}+A^{*}A \Vert .
$$
\end{theorem}
\noindent {\it Proof.}
For each $ t\in [0,1] $ and $ x\in S_{1}(H) $ using Lemma \ref{lem 3.4} we have 
\begin{eqnarray}
& & \vert \langle( \varphi (t)A+\psi(t)A^{*})x,x\rangle \vert^{2} \nonumber \\
& \leq &  \left( \vert \varphi \vert (t) \vert \langle Ax,x\rangle \vert + \vert \psi \vert (t) \vert \langle A^{*}x,x\rangle \vert \right) ^{2} \nonumber \\
& = & \vert \varphi \vert^{2} (t) \vert \langle Ax,x\rangle \vert^{2}+ 2\vert \varphi \psi \vert (t) \vert \langle Ax,x\rangle \langle x, A^{*}x\rangle \vert + \vert \psi \vert^{2} (t) \vert \langle Ax,x\rangle \vert^{2} \nonumber \\
& = & \left( \vert \varphi \vert^{2} + \vert \psi \vert^{2}\right)(t) \vert \langle Ax,x\rangle \vert^{2}+ \vert \varphi \psi \vert (t) \vert \langle Ax, A^{*}x\rangle \vert + \vert \varphi \psi \vert (t) \Vert Ax \Vert \Vert A^{*} x\Vert \nonumber \\
& \leq & \left( \vert \varphi \vert^{2} + \vert \psi \vert^{2}\right)(t) \vert \langle Ax,x\rangle \vert^{2}+ \vert \varphi \psi \vert (t) \vert \langle A^{2}x,x\rangle \vert + \frac{1}{2}\vert \varphi \psi \vert (t) \left(  \Vert Ax \Vert^{2} +  \Vert A^{*}x \Vert^{2} \right)  \nonumber \\
& \leq & \left( \vert \varphi \vert^{2}+ \vert \psi \vert^{2} \right)(t)\omega^{2}(A)+\vert \varphi \psi \vert (t) \omega (A^{2})+ \frac{1}{2}\vert \varphi \psi \vert (t) \Vert AA^{*}+A^{*}A \Vert . \nonumber
\end{eqnarray}
From the last estimates and definition of the numerical radius validity of theorem is clear.
\begin{corollary}
\label{cor3.9}
If we take $ \varphi (t)=1-2t, \ \psi (t)=1, \ 0\leq t \leq 1 $ in Theorem \ref{thm3.11}, we get the inequality proved in [26, Th. 2.9], which says that
$$
\omega_{t}^{2}(T)\leq (2-4t+4t^{2})\omega^{2}(T)+(1-2t)\omega(T^{2})+\frac{1}{2}(1-2t)\Vert T^{*}T+TT^{*}\Vert.
$$
So, the inequality obtained in Theorem 3.8 generalizes [26, Th. 2.9].
\end{corollary}
\begin{remark}
Let  $ A\in \mathbb{B}(H) $ and $ t\in [0,1], $ then from Theorem \ref{thm3.11} we have
\begin{eqnarray}
\omega^{2}_{t}(\varphi ,\psi;A) & \leq & \left( \vert \varphi \vert^{2}+ \vert \psi \vert^{2} \right)(t)\omega^{2}(A)+\vert \varphi \psi \vert (t) \omega (A^{2})+ \frac{1}{2}\vert \varphi \psi \vert (t) \Vert AA^{*}+A^{*}A \Vert  \nonumber \\
& \leq & \left( \vert \varphi \vert^{2}+ \vert \psi \vert^{2} \right)(t)\Vert A \Vert^{2}+\vert \varphi \psi \vert (t) \Vert A \Vert^{2} + \frac{1}{2}\vert \varphi \psi \vert (t) 2\Vert A \Vert^{2}  \nonumber \\
& = & \left( \vert \varphi \vert + \vert \psi \vert \right) ^{2}\Vert A \Vert^{2}. \nonumber
\end{eqnarray}
Clearly, for $ \varphi (t)=0, \ \psi (t)=1, \ 0\leq t \leq 1 $ the upper bound improves the second inequality in (\ref{equ1}). So, the inequality obtained in Theorem \ref{thm3.11} generalizes the second ineguality in (\ref{equ1}).
\end{remark}
\noindent {\it Example 2.} In the complex Hilbert space $ L^{2}(-1,1) $, consider the following Skew-symmetric Volterra integration operator
$$
A:L^{2}(-1,1)\rightarrow L^{2}(-1,1), \ Af(x)=\int\limits_{-x}^{x} f(t)dt, \ f\in L^{2}(-1,1).
$$
It is known that $ A $ is a nilpotent operator with index, $ \Vert A \Vert=\frac{4}{\pi}, \ \omega (A)=\frac{2}{\pi} $ (see \cite{Halmos}, \cite{Karaev}).

On the other hand, in work \cite{Kittaneh} Kittaneh proved that if $ A^{2}=0, $ then $ \Vert A \Vert^{2}= \Vert  A^{*}A+AA*\Vert .$

Using these results by Theorem \ref{thm 3.7}, \ref{thm3.9}, and \ref{thm3.11}, we have
$$
\omega_{t}^{2}(\varphi,\psi;A)\leq \frac{4}{\pi^{2}}\left( \vert \varphi \vert^{2}(t)+2\vert \psi \vert^{2}(t)+2\vert \varphi \psi \vert(t)\right) , 
$$
$$
\omega_{t}^{2}(\varphi,\psi;A)\leq \frac{8}{\pi^{2}}\left( \vert \varphi \vert^{2}+\vert \psi \vert^{2}\right)(t) , 
$$
and
$$
\omega_{t}^{2}(\varphi,\psi;A)\leq \frac{4}{\pi^{2}}\left( \vert \varphi \vert^{2}+\vert \psi \vert^{2}\right)(t)+ \frac{8}{\pi^{2}} \vert \varphi \psi \vert(t) = \frac{4}{\pi^{2}}\left( \vert \varphi \vert+\vert \psi \vert\right)^{2}(t),
$$
respectively.

\section{Weighted numerical index in complex Hilbert spaces}
In this section, the lower and upper evaluations for the weighted numerical index in Hilbert spaces are given.

\begin{theorem}
\label{thm4.1}
Let $ A\in \mathbb{B}(H) $ and $ t\in [0,1],$ then the weighted numerical radius function and the weighted numerical index hold the following inequalities
$$
\inf\left\lbrace \vert \varphi + \psi \vert (t), \vert \varphi - \psi \vert (t)  \right\rbrace \omega (A) \leq \omega_{t} (\varphi, \psi; A) \leq \sup\left\lbrace \vert \varphi + \psi \vert (t), \vert \varphi - \psi \vert (t)  \right\rbrace \omega(A)
$$
and 
$$
\frac{1}{2}\inf\left\lbrace \vert \varphi + \psi \vert (t), \vert \varphi - \psi \vert (t)  \right\rbrace \leq n _{t}(\varphi, \psi; H) \leq \frac{1}{2}\sup\left\lbrace \vert \varphi + \psi \vert (t), \vert \varphi - \psi \vert (t)  \right\rbrace ,
$$
respectively.
\end{theorem}
\noindent {\it Proof.}
For the $ A\in \mathbb{B}(H) $ and $ \ x\in S_{1}(H) ,$ we have 
$$
\varphi A+ \psi A^{*}= (\varphi + \psi )ReA +i(\varphi - \psi)Im A
$$
and 
\begin{eqnarray}
 \inf\left\lbrace \vert \varphi + \psi \vert (t), \vert \varphi - \psi \vert (t)  \right\rbrace \vert \langle Ax,x\rangle \vert 
& = & \inf\left\lbrace \vert \varphi + \psi \vert (t), \vert \varphi - \psi \vert (t)  \right\rbrace \sqrt{\langle Re Ax,x\rangle^{2}+\langle ImAx,x\rangle^{2}} \nonumber \\
& \leq & \sqrt{\langle\vert \varphi + \psi \vert (t)Re Ax,x\rangle^{2}+\langle\vert \varphi - \psi \vert (t)ImAx,x\rangle^{2}} \nonumber .
\end{eqnarray}
Then,
$$
\inf\left\lbrace \vert \varphi + \psi \vert (t), \vert \varphi - \psi \vert (t)  \right\rbrace  \omega (A) \leq \omega_{t} (\varphi, \psi; A).
$$
On the other hand, for $ x\in S_{1}(H) $ since 
\begin{eqnarray}
\vert \langle(\varphi A+\psi A^{*})x,x\rangle \vert & = & \vert \langle(\varphi + \psi )ReA+i(\varphi -\psi)ImA)x,x\rangle  \vert \nonumber \\
& = & \sqrt{\langle( \varphi + \psi ) (t)Re Ax,x\rangle^{2}+\langle( \varphi - \psi ) (t)ImAx,x\rangle^{2}} \nonumber \\
& \leq & \sup\left\lbrace \vert \varphi + \psi \vert (t), \vert \varphi - \psi \vert (t)  \right\rbrace \sqrt{\langle Re Ax,x\rangle^{2}+\langle ImAx,x\rangle^{2}} \nonumber \\
& = & \sup\left\lbrace \vert \varphi + \psi \vert (t), \vert \varphi - \psi \vert (t)  \right\rbrace  \vert \langle Ax,x\rangle \vert , \nonumber
\end{eqnarray}
then 
$$
\omega_{t} (\varphi, \psi; A) \leq \sup\left\lbrace \vert \varphi + \psi \vert (t), \vert \varphi - \psi \vert (t)  \right\rbrace \omega(A).
$$
From the first claim of this theorem, the definition of the numerical index, and the result $ n(H)=\frac{1}{2} $ \cite{Bonsall1}, \cite{Kadets} the second claim of this theorem is obtained.
\begin{corollary}
\label{cor 4.2}
Using the inequalities the first claim of Theorem \ref{thm4.1} and (\ref{equ2}) it can be proved that for any $ A\in \mathbb{B}(H) $ and $ t\in [0,1] $ the following inequalities hold
$$
\frac{1}{4}\left( \inf\left\lbrace \vert \varphi + \psi \vert (t), \vert \varphi - \psi \vert (t)  \right\rbrace \right) ^{2}\Vert A^{*}A+AA^{*}\Vert \leq \omega_{t}^{2}(\varphi,\psi;A)\leq \frac{1}{2}\left( \sup\left\lbrace \vert \varphi + \psi \vert (t), \vert \varphi - \psi \vert (t)  \right\rbrace \right) ^{2}\Vert A^{*}A+AA^{*}\Vert.
$$
\end{corollary}
\begin{corollary}
If we take $ \varphi (t)=1, \ \psi (t)=0, \ 0\leq t \leq 1 $ in Corollary \ref{cor 4.2}, we get the inequality (\ref{equ2}). So, the inequality obtained in Corollary \ref{cor 4.2} generalizes (\ref{equ2}).
\end{corollary}
\begin{corollary}
Let $ A\in \mathbb{B}(H), $ then
\begin{eqnarray}
\int\limits_{0}^{1} \inf\left\lbrace \vert \varphi + \psi \vert (t), \vert \varphi - \psi \vert (t)  \right\rbrace dt \omega (A)
 \leq  \int\limits_{0}^{1}  \omega_{t}(\varphi, \psi; A)dt
 \leq  \int\limits_{0}^{1} \sup\left\lbrace \vert \varphi + \psi \vert (t), \vert \varphi - \psi \vert (t)  \right\rbrace dt \omega (A). \nonumber
\end{eqnarray}
\end{corollary}
\begin{corollary}
If we take $ \varphi (t)=1, \ \psi(t)=0, \ t\in [0,1], $ we get
\begin{eqnarray}
\inf\left\lbrace \vert \varphi + \psi \vert (t), \vert \varphi - \psi \vert (t)  \right\rbrace =1 \ \ \text{and} \ \ \ \sup\left\lbrace \vert \varphi + \psi \vert (t), \vert \varphi - \psi \vert (t)  \right\rbrace=1. \nonumber
\end{eqnarray}
By the first claim of Theorem 4.1 we obtain $ n(1,0;H)=\frac{1}{2}, \ t\in [0,1], $ which founds in  \cite{Bonsall1} and \cite{Kadets}.
\end{corollary}
\begin{corollary}
If we take $ \varphi (t)=1, \ \psi (t)= 1-2t, \ t\in [0,1], $ we get
$$
\vert \varphi + \psi  \vert (t)=2-2t, \ \vert \varphi - \psi  \vert (t)=2t, \ t\in [0,1],
$$
then 
$$
\alpha (t)= \inf\left\lbrace \vert \varphi + \psi \vert (t), \vert \varphi - \psi \vert (t)  \right\rbrace= \inf\left\lbrace 2-2t,2t \right\rbrace =
\begin{cases}
2t, & 0\leq t\leq \frac{1}{2}, \\
2-2t, & \frac{1}{2}\leq t \leq 1
\end{cases}
$$
and 
$$
\lambda (t)= \sup\left\lbrace \vert \varphi + \psi \vert (t), \vert \varphi - \psi \vert (t)  \right\rbrace= \sup\left\lbrace 2-2t,2t \right\rbrace =
\begin{cases}
2-2t, & 0\leq t\leq \frac{1}{2}, \\
2t, & \frac{1}{2}\leq t \leq 1
\end{cases}.
$$
Thus, by the Theorem \ref{thm4.1}
$$
\frac{1}{2}\alpha (t) \leq n(1,1-2t;H)\leq \frac{1}{2}\lambda (t) , \ t\in [0,1].
$$
In special case, when $ t=\frac{1}{2},  $ from the above relations we have
$$
n_{t}(H)=n (1,0;H) = \frac{1}{2}, \ t\in [0,1]
$$
which founds in \cite{Bonsall1} and \cite{Kadets}.
\end{corollary}
\begin{corollary}
If we take $ \varphi (t)=\sin\left( \frac{1}{4}t\right) , \ \psi(t)=\cos\left( \frac{1}{4}t\right), \ t\in [0,\pi], $
then
$$
 \inf\left\lbrace \vert \varphi + \psi \vert (t), \vert \varphi - \psi \vert (t)  \right\rbrace= \sin\left( \frac{1}{4}t\right) \ \ \text{and} \ \
 \sup\left\lbrace \vert \varphi + \psi \vert (t), \vert \varphi - \psi \vert (t)  \right\rbrace= \cos\left( \frac{1}{4}t\right), \ t\in [0,\pi ].
$$
So, by the Theorem \ref{thm4.1} we have 
$$
\frac{1}{2}\sin\left( \frac{1}{4}t\right) \leq n(\varphi,\psi;H)\leq \frac{1}{2}\cos\left( \frac{1}{4}t\right), \ t\in [0,\pi ].
$$
In special value of $ t=\pi $ it is obtained that 
$$
n(\varphi(\pi),\psi(\pi);H)=n(\frac{\sqrt{2}}{2},\frac{\sqrt{2}}{2};H)=2^{-\frac{3}{2}}
$$
\end{corollary}

Now, we give one approximation for the weighted numerical index.
\begin{theorem}
\label{thm4.8}
For any $ A\in \mathbb{B}(H) $ the followings are true
\begin{enumerate}[label=(\arabic*)]
\item If $ \varphi_{1}, \ \varphi_{2}, \ \psi_{1}, \ \psi_{2}\in C[0,1], $ then
$
\vert n_{t}(\varphi_{1},\psi_{1};H)-n_{t}(\varphi_{2},\psi_{2};H)\vert \leq \vert \varphi_{1}-\varphi_{2} \vert (t)+\vert \psi_{1}-\psi_{2} \vert (t), \ t\in [0,1],
$
\item If $ \varphi, \ \psi, \varphi_{n}, \ \psi_{n}\in C[0,1] $ and sequences $ (\varphi_{n}) $ and $ (\psi_{n}) $ converge to the functions $ \varphi $ and $ \psi $, respectively, then
$$
n_{t}(\varphi,\psi;H)=\lim\limits_{n\rightarrow \infty} n_{t}(\varphi_{n},\psi_{n};H).
$$
\end{enumerate}
\end{theorem}
\noindent {\it Proof.}
(1) For $ x\in H, \ \Vert x \Vert =1 $, we have 
\begin{eqnarray}
\vert \langle \varphi_{1}(t)A+\psi_{1}(t)A^{*})x,x\rangle\vert  \leq  \vert \langle (\varphi_{1}-\varphi_{2})(t)A+(\psi_{1}-\psi_{2})(t)A^{*})x,x\rangle\vert 
 +  \vert \langle \varphi_{2}(t)A+\psi_{2}(t)A^{*})x,x\rangle\vert . \nonumber
\end{eqnarray}
Hence, for $ \Vert A \Vert =1 $ we get
$$
\omega (\varphi_{1},\psi_{1};A)\leq \sup\left\lbrace \vert  \varphi_{1}-\varphi_{2} \vert (t), \vert  \psi_{1}-\psi_{2} \vert (t)\right\rbrace \Vert A \Vert  +\omega (\varphi_{2},\psi_{2};A).
$$
Consequently, it is obtained that 
$$
n(\varphi_{1},\psi_{1};H)\leq \sup \left\lbrace \vert\varphi_{1}-\varphi_{2} \vert (t),\vert\psi_{1}-\psi_{2} \vert (t)  \right\rbrace  + n(\varphi_{2},\psi_{2};H).
$$
Similarly, it can be shown that 
$$
n(\varphi_{2},\psi_{2};H)\leq \sup \left\lbrace \vert\varphi_{1}-\varphi_{2} \vert (t),\vert\psi_{1}-\psi_{2} \vert (t)  \right\rbrace  + n(\varphi_{1},\psi_{1};H).
$$
So, we have 
$$
\vert n(\varphi_{1},\psi_{1};H)-n(\varphi_{2},\psi_{2};H)\vert \leq \vert \varphi_{1}-\varphi_{2} \vert (t)+\vert \psi_{1}-\psi_{2} \vert (t), \ t\in [0,1].
$$
(2) From the first claim of this theorem we have 
$$
\vert n(\varphi_{n},\psi_{n};H)-n(\varphi,\psi;H)\vert \leq \vert \varphi_{n}-\varphi \vert (t)+\vert \psi_{n}-\psi \vert (t), \ t\in [0,1], \ n\geq 1.
$$
Since $ \varphi_{n}\rightarrow \varphi $ and $ \psi_{n}\rightarrow \psi  $ as $ n\rightarrow \infty , $ then from the last relation the validity of second claim it is clear.

\begin{corollary}
If $ \varphi_{n}, \psi_{n}\in C[0,1], n\geq 1 $ and $ \varphi_{n}\rightarrow 1, \psi_{n}\rightarrow 0 $ as $ n\rightarrow \infty, $ then from the claim (2) of Theorem \ref{thm4.8} we have 
$$
\lim\limits_{n\rightarrow \infty} n_{t}(\varphi_{n},\psi_{n};H)=\frac{1}{2}.
$$
\end{corollary}

\end{document}